\font\Ch=msbm10  
\font\ch=msbm8  
\def \QQ {\hbox{\Ch Q}}
\def \ZZ {\hbox{\Ch Z}}
\def \zz {\hbox{\ch Z}}

\def\la{{\lambda}}
\def\htq #1{{\widetilde{#1}^q}}
\def\htt #1{{\widetilde{#1}^t}}
\overfullrule=0pt
\baselineskip 14pt
\settabs 8 \columns
\+ \hfill&\hfill \hfill & \hfill \hfill   & \hfill \hfill  & \hfill \hfill   &
 & \hfill \hfill\cr
\parindent=.5truein
\hfuzz=3.44182pt
\hsize 6truein
\font\ita=cmsl12 
\font\small=cmr6
\font\title=cmbx10 scaled\magstep2
\font\normal=cmr12 
\font\ninerm=cmr9 
\font\sbol=cmbx9
\font\bol=cmbx12

\def \con {\subseteq}

\def\sig{\sigma}

\def \-> {\rightarrow}

\def\OM {\Omega}

\def\la {\lambda}
\def\La {\Lambda}
\def \RA {\rightarrow}

\def\xon {x_1,x_2,\ldots ,x_n}

\def \sas {\vskip .06truein}
\def\sa{{\vskip .125truein}}

\def\sap{{\vskip .25truein}}

\def \eee {\epsilon}

\def\con {\subseteq}
\def \ses {\enskip = \enskip}
\def \sps {\enskip + \enskip}

\def \scs {\ssp , \ssp}
\def \ess {\enskip}
\def \ssp {\hskip .25em}

\def \part {\vdash}

\normal

\vsize=8truein
\sap
\def\today{\ifcase\month\or
January\or February\or March\or April\or may\or June\or
July\or August\or September\or October\or November\or
December\fi
\space\number\day, \number\year}

\headline={   
\small$q,t$-Kostka Polynomiality Revisited $\ess\ess\ess\ess\ess$ 
\hfill$\ess\ess\ess$\today$\ess\ess\ess\ess\ess\ess$ \folio }
 \footline={\hfil}

\def \OM {{\Omega}}
\def\OH {{G}}

\def \RA {{ \rightarrow }}

\def \TK {{\tilde K}}
\def \TH {{\tilde H}}

\def \xon {x_1,\ldots ,x_n}

\def \TK {{\tilde K}}
\def \TH {{\tilde H}}
\def \xon {x_1,\ldots ,x_n}

\font\small=cmr6
\def \scs {\ssp , \ssp}
\def \ess {\enskip}
\def \ssp {\hskip .25em}

\def \part {\vdash}
\font\title=cmbx10 scaled\magstep2
\font\normal=cmr10 
\def\today{\ifcase\month\or
January\or February\or March\or April\or may\or June\or
July\or August\or September\or October\or November\or
December\fi
\space\number\day, \number\year}
\centerline{\bol Polynomiality of the $\bf q,t$-Kostka Revisited}
\vskip .1in
\centerline{  by}
\vskip .1in
\centerline{\bf Adriano M. Garsia \footnote{$\dag$}{\rm Supported by NSF}}
\centerline{\it  University of California San Diego  }
\vskip .1in
\centerline{  and}
\vskip .1in
\centerline{\bf Mike Zabrocki  }
\centerline{\it{ Universit\'e du Qu\'ebec \`a Montr\'eal}}

\narrower{
\noindent{\sbol Abstract.}

{\ninerm
Let $K(q,t)= \|K_{\la\mu}(q,t)\|_{\la,\mu}$ be the Macdonald $q,t\,$-Kostka matrix and $K(t)=K(0,t)$
be the matrix of the Kostka-Foulkes polynomials $K_{\la\mu}(t)$. 
In this paper we present a new proof of the polynomiality of the $q,t$-Kostka coefficients
that is both short and elementary. More precisely, we derive that $K(q,t)$ has entries in $ \ZZ[q,t]$
directly from the fact that the matrix $K(t)^{-1}$ has entries in 
$\ZZ[t]$.  The proof uses only identities that can be found in the original
paper [7] of Macdonald.  }
}
\sa

\hsize=6.8truein
\noindent{\bol Introduction}
\sas

The polynomiality problem for the  $q,t$-Kostka coefficients [11], was posed by Macdonald in 
the fall 1988 meeting of the Lotharingian seminar. It remained open for quite a few years, when suddenly
in 1996, several proofs of varied difficulty appeared in a period of only a few months.
At the present there are three basically  different approaches to proving the polynomiality
of the $q,t$-Kostka coefficients: 
\sas

\item{(1)} Via plethystic formulas (Garsia-Tesler [4],  Garsia-Remmel [3]).
\sas

\item{(2)}   Via vanishing properties (Sahi [13] and Knop [7],[8]).
\sas

\item{(3)}  Via Rodriguez formulas (Lapointe-Vinet [10], Kirillov-Noumi [6])
\sas

Each of these approaches has its own special advantages. The plethystic approach
led to very efficient algorithms for computing these coefficients and ultimately
produced some remarkably simple explicit formulas [2]. The vanishing properties 
approach led to the discovery some  basic non-symmetric variants of the 
Macdonald  polynomials with remarkable combinatorial implications that still remain
to be fully explored. The approach via Rodriguez formulas stems from a pioneering
paper of Lapointe-Vinet [9] on Jack-Polynomials. Although originally it was
based on   deep affine Hecke algebra identities, eventually the idea led to
some of the most elementary proofs of the polynomiality result (see [10] and
[6]). In particular it produced a family of symmetric  function operators
$\{B_k^{q,t}\}_{k=1,2,\ldots}$, which permitted the construction of the 
Macdonald ``integral forms''
$J_\mu(x;q,t)$, one part at the time, starting from $1$, according to  an
identity of the form
$$
J_{\mu'}(x;q,t)\ses B_{\mu_1}^{q,t} B_{\mu_2}^{q,t} \cdots B_{\mu_k}^{q,t}\ssp
{\bf 1}\ess .
$$ 
Our main contribution here is a remarkably simple argument which shows that families of operators
$B_k^{q,t}$ yielding such a formula, may obtained by $q$-twisting in a minor
way {\bf any} sequence  of operators $\{B_k^{t}\}_{k=1,2,\ldots}$ which yields
the analogous formula 
$$
Q_{\mu'}(x;t)\ses B_{\mu_1}^{ t} B_{\mu_2}^{ t} \cdots B_{\mu_k}^{ t} \ssp {\bf
1}\ess .
$$ 
for the Hall-Littlewood polynomial $Q_\mu(x;q,t)$.
 As a byproduct we obtain that the polynomiality
of the $q,t$-Kostka is an immediate consequence of the polynomiality of the the Kostka-Foulkes coefficients.
What is surprising is that this fact was missed for so  many years by researchers in this area.
What might be even more surprising is that we obtain a remarkably general result by further 
simplifying some of the arguments used in [10] and [6].  To give a more precise
description of our results   we need some notation.  
\sas

We shall deal with identities in the  algebra $\La$ of symmetric functions in a finite or infinite
alphabet $X=\{x_1,x_2,x_3,\ldots \}$ with
coefficients  in the field of rational functions $\QQ(q,t)$. We also denote by
$\La_{\zz[q,t\, ]}$ the algebra of symmetric functions  in $X $
with coefficients  in $\ZZ[q,t\, ]$
 
An essential notational tool in our presentation is the notion of ``{\ita plethystic substitution}''  and we need to recall its definition.
Briefly, if 
$E=E(t_1,t_2,t_3,\ldots )$ is a given formal series
in the variables $t_1,t_2,t_3,\ldots $ (which may include the parameters $q,t$)
and $f\in \Lambda$ has
been expressed in terms of the power basis in the form
$$
F\ses Q(p_1,p_2,p_3,\ldots \, )
$$
then the ``{\ita plethystic substitution}'' of $E$ in $F$, denoted $F[E]$, is  simply defined by
setting
$$
F[E]\ses Q(p_1,p_2,p_3,\ldots \, )\, \Big|_{p_k \RA E(t_1^k,t_2^k,t_3^k,\ldots
)}\ess .
\eqno  {\rm I}.1
$$
This operation is easily programmed in any symbolic manipulation software which 
includes a symmetric function package.
We shall  $\,$ adopt the convention that inside the plethystic brackets
$[\,\,]$, $X$ and $X_n$ respectively stand for $x_1+x_2+x_3+\cdots $ and 
$x_1+x_2+\cdots + x_n$. A similar notation is adopted when we work with any
other alphabet. Note that if $P\in\Lambda$ then $P[X_n]$ simply means
$P(x_1,x_2,\ldots ,x_n)$. Plethystic substitution need not be restricted to
symmetric polynomials. In fact, the substitution in I.1 makes sense even  when
$Q$ is a formal power series.  In this vein, we shall systematically use the
symbol $\OM[X]$ to represent the symmetric function series
$$
\OM[X]\ses \prod_{i} {1\over 1-x_i}\ses \exp\Big(\sum_{k\geq 1}{p_k[X]\over k}\Big)\ess .
\eqno {\rm I}.2
$$
This notation is particularly convenient since for any two alphabets $X,Y$
from I.1 and I.2 we easily derive that
$$
\OM\big[ X+Y\big]\ses \OM\big[ X ]\OM\big[  Y\big] 
\scs\ess\ess\ess 
\OM\big[ X-Y\big]\ses \OM\big[ X ]/\OM\big[  Y\big] 
\eqno {\rm I}.3
$$
In particular we see that the Macdonald kernel $\OM_{q,t}(x_1,x_2,\ldots ,x_n;y_1,y_2,\ldots ,y_k)$ 
(see [12] (2.5) p.~309) may be expressed in the compact form.
$$
\OM_{q,t}(x_1,x_2,\ldots ,x_n;y_1,y_2,\ldots ,y_k)\ses \OM\big[ X_nY_k\textstyle{1-t\over 1-q}\big]
\eqno {\rm I}.4
$$

We represent partitions here by their French Ferrers diagram, that is
with rows decreasing from bottom to top. To express that a certain lattice cell $s$ belongs to the 
Ferrers diagram of a partition $\la$ we shall simply write $s\in \la$. More
generally we shall always identify partitions with their corresponding
diagrams. The number of rows of the Ferrers diagram of a partition
$\mu$ will be called the {\ita length} of $\mu$ and denoted $l(\mu)$. We
shall often make use of the operation of prepending a column of length $k$ to
the diagram of a partition $\mu$ of length $\leq k$, for convenience we shall
denote the resulting partition by the symbol $\mu+1^k$.
\sas

Let us recall that the {\ita integral forms} $J_\la(x;q,t)$ of the
Macdonald   polynomials $P_\la(x;q,t)$ and $Q_\la(x;q,t)$  are defined (see
[12] p.~352) by setting
$$
J_\la(x;q,t)=h_\la(q,t)P_\la(x;q,t)=h_\la'(q,t)Q_\la(x;q,t)
\eqno {\rm I}.5
$$
with
$$
h_\la(q,t)\ses \prod_{s\in\la} \big(1-q^{a_\la(s)}t^{l_\mu(s)+1}\big)
\scs
\ess\ess\ess\ess
h_\la'(q,t)\ses \prod_{s\in\la} \big(1- q^{a_\la(s)+1} t^{l_\mu(s) }\big)
\eqno {\rm I}.6
$$
where, for a cell $s\in\la$, $a_\la(s)$ and $l_\la(s)$ represent the {\ita
arm} and {\ita leg} of $s$ in $\la$, that is  the number of cells of
$\la$ that are respectively strictly {\it  EAST} and {\it  NORTH} of $s$. 
Recall that Garsia-Haiman [1] introduce the
modified versions $H_\mu[X;q,t\, ]$ of the integral forms by setting
$$
H_\mu[X;q,t\, ]\ses J_\mu[\textstyle{X\over  1-t};q,t\, ]
\eqno {\rm I}.7
$$
These polynomials offer direct access to the Macdonald $q,t$-Kostka
coefficients because their Schur function
expansion reduces to
$$
H_\mu[X;q,t\, ]\ses\sum_{\la} S_\la[X] K_{\la\mu}(q,t)
$$
This follows immediately from the definition of $K_{\la\mu}(q,t)$, given in [12] 
((8.11) p.~354), as the coefficients appearing in the expansion
$$
J_\mu[X;q,t\, ]\ses\sum_{\la} S_\la[X(1-t)]\,  K_{\la\mu}(q,t)
\eqno {\rm I}.8
$$ 
It is also shown in [12] ((8.12) p.~354) that we have
$$
J_\mu[X;0,t\, ]\ses Q_\mu[X;t\, ] = \sum_{\la} S_\la[X(1-t)]\,  K_{\la\mu}(t)
\eqno {\rm I}.9
$$
where $Q_\la[X;t\, ]$ is the Hall-Littlewood polynomial and $K_{\la\mu}(t)$ denotes the corresponding
{\ita Kostka-Foulkes} coefficient. We shall also set
$$
H_\mu[X;t\, ]\ses Q_\la[{\textstyle{X \over 1-t}};t\,] \ses \sum_{\la}\, S_\la[X] K_{\la\mu}(t)\,\, . 
\eqno {\rm I}.10
$$

It develops that the polynomiality of the $K_{\la\mu}(q,t)$ is an immediate consequence of the
following general result first proved in [14].
\sa

\noindent{\bol Theorem I.1}

{\ita
 For any linear operator $V$ acting on $\Lambda$ and $P\in \Lambda$   set
$$
\htq{V}P[X]\ses V^YP[qX+(1-q)Y]\, \big|_{Y=X}
\eqno {\rm I}.11
$$
where $V^Y$ is simply $V$ acting on polynomials in the $Y$ variables. This given,  
if $\OH_k=\OH_k(X,t)$ is any linear operator on $\Lambda$ with the property that
$$
\OH_k H_\mu[X;t\, ]\ses H_{\mu+1^k} [X;t\, ]
\ess\ess\ess\ess \hbox{for all $\mu$ of length $\leq k$} 
\eqno {\rm I}.12
$$
then $\htq{\OH_k}$ has the property
$$
\htq{\OH_k} H_\mu[X;q,t\, ]\ses H_{\mu+1^k} [X;q,t\, ]
\ess\ess\ess\ess \hbox{for all $\mu$ of length $\leq k$}\ess . 
\eqno {\rm I}.13
$$
In particular, the modified Macdonald polynomial $H_\mu[X;q,t\,]$ may be obtained from the 
``Rodriguez'' formula:
$$
H_\mu[X;q,t\,]\ses \htq{\OH_{\mu_1'}}\htq{\OH_{\mu_2'}}\cdots \htq{\OH_{\mu_h'}}\,
{\bf 1}
\eqno {\rm I}.14
$$
where $\mu'=(\mu_1',\mu_2',\ldots ,\mu_h')$ denotes the conjugate of $\mu$.
}
\sa

Our main contribution here is a simple, direct and elementary proof of this result 
which only uses identities given in the  original paper of Macdonald.
\sas

Note that since the Kostka-Foulkes matrix $K(t)=\|K_{\la\mu}(t)\|_{\la\mu}$
is unitriangular it follows that its inverse $H(t)=K(t)^{-1}$ has entries in
$\ZZ[t]$. This implies that the ``trivial'' operator $T\OH_k=T\OH_k(X;t)$ defined
by setting for the $\{H_\mu[X;t]\}_\mu$ basis  
$$
T\OH_k H_\mu[X;t\, ]\ses 
\cases{
H_{\mu+1^k} [X;t\, ]& if  $\ess \rm l(\mu)\leq k$\cr\cr 
0 & $\rm otherwise$
}
\eqno {\rm I}.15
$$
acts integrally on the Schur basis.  
This given, we see that the desired result
$$
K_{\la\mu}(q,t)\in \ZZ[q,t]
\eqno {\rm I}.16
$$
is an immediate consequence of I.14 with $\OH=T\OH$.
\sa

This paper is divided in two sections. In the first section we prove Theorem I.1,
and in the second section we give a number of applications, including the explicit derivation of
the action of a variant of the operator $T\OH$ on the monomial basis.
\sa

\noindent{\bol 1. Rodriguez operators for the Integral Forms}
\sas

We shall start by proving a result analogous to Theorem I.1 
for the Macdonald integral forms. To this end we shall need a number of 
auxiliary results. 
\sap

\noindent{\bol Proposition 1.1}

{\ita If $V(X;t)$ is a linear operator on symmetric functions in $X$ that
depends only on $t$ and we set
$$
\htq{V}(X;t)  P[X]\ses
V(Z;t)^ZP\big[qX+(1-q)Z\big]\Big|_{Z=X}
\eqno 1.1
$$
then  
$$
{\htq{V}(X;t)\OM\left[XY{\textstyle{1-t\over 1-q}}\right]\over
\OM\left[XY{\textstyle{1-t\over 1-q}}\right]}
\ses
  {V (X;t)^X\OM\left[   X Y{\textstyle(1-t)}\right] \over  \OM\left[   X Y{\textstyle(1-t )}\right] }
\eqno 1.2
$$
in particular this ratio is independent of $q$.}

\noindent {\bol Proof}

We have
$$
\eqalign{
\htq{V}(X;t)\OM\left[XY{\textstyle{1-t\over 1-q}}\right]
&\ses 
V (Z;t)^Z\OM\left[\big(qX+(1-q)Z\big)Y{\textstyle{1-t\over
1-q}}\right]\Big|_{Z=X}\cr &\ses 
\OM\left[  qX Y{\textstyle{1-t\over 1-q}}\right] V (X;t)^X\OM\left[\big (1-q)X
Y{\textstyle{1-t\over 1-q}}\right]\cr &\ses 
\OM\left[  X Y{\textstyle{1-t\over 1-q}}\right] \OM\left[  (q-1)X Y{\textstyle{1-t\over 1-q}}\right]
  V (X;t)^X\OM\left[   X Y{\textstyle(1-t
)}\right] \cr 
&\ses 
\OM\left[  X Y{\textstyle{1-t\over 1-q}}\right] 
\OM\left[   X Y{\textstyle(t-1 )}\right]
  V (X;t)^X\OM\left[   X Y{\textstyle(1-t
)}\right] \cr}
$$
This proves 1.2.
\sas

Recall that for an alphabet $Y_k=y_1+y_2+\cdots +y_k$ and an interval $I\con [1,k]=\{1,2,\ldots k\}$ 
Macdonald sets  (see (3.5) [12] p.~315).
$$
A_I[Y_k;t]\ses t^{|I|\choose 2}\prod_{i\in I\, \atop j\in\,  [1,k]-I } {t
y_i-y_j\over  y_i-y_j }
$$
Note also that for any symmetric polynomial $P[Y_k]$, the operator $T_I^q$ which replaces $y_i$ by
$qy_i$ for every $i\in I$ may be written in the form
$$
T_I^qP[Y_k]\ses P[Y_k+(q-1)Y_I]
\eqno 1.3
$$
where 
$$
Y_I\ses \sum_{i\in I}y_i\ess .
$$
This given, we have the following identity due to Kirillov-Noumi [6].
\sas

\noindent{\bol Proposition 1.2} 

{\ita Let $M^{(Y_k)}(u)$ denote the Macdonald operator acting on polynomials
in the alphabet $Y_k$. That is  
$$
M^{(Y_k)}(u)\ses \sum_{r=0}^k  u^r\sum_{I\con [1,k] \atop
|I|=r}  A_I[Y_k;t] T_I^q
$$ 
then for $k\leq n$}
$$
M^{(Y_k)}(-1)\OM\left[  X_n Y_k{\textstyle{1-t\over 1-q}}\right] \ses 
y_1y_2\cdots y_k\sum_{l(\mu)\leq k }J_{\mu+1^k}[X_n;q,t]\,{ P_\mu[Y_k;q,t]\over h_\mu'(q,t)}
\eqno 1.4
$$
\noindent{\bol Proof}

We start with the Macdonald {\ita Cauchy}-identity from equation (4.13) in
[12]
$$
\OM\big[X_nY_k{\textstyle{1-t\over 1-q}}\big]\ses \sum_{l(\la)\leq k} Q_{\la}[X_n ;q,t\, ]
P_\la [Y_k;q,t\, ]
$$
the summation being only over partitions $\la$ of length $\leq k$ because $P_\la$ vanishes
when evaluated on an alphabet whose cardinality is smaller than the length of $\la$. 
Applying $M^{(Y_k)}(u)$ to both sides of this identity and using Theorem (4.15) of
[12] (p.~324) we get that
$$
M^{(Y_k)}(u)\OM\big[X_nY_k{\textstyle{1-t\over 1-q}}\big]\ses \sum_{l(\la)\leq k} Q_{\la}[X_n ;q,t\, ]
\Big(\prod_{i=1}^k (1+u\, t^{k-i}q^{\la_i})\Big)P_\la [Y_k;q,t\, ]\, ,
$$
Now note that if $l(\la)<k$ then the term corresponding to $\la$ in this sum will vanish if we set $u=-1$.
We thus obtain
$$
M^{(Y_k)}(-1)\OM\big[X_nY_k{\textstyle{1-t\over 1-q}}\big]\ses \sum_{l(\la)= k} Q_{\la}[X_n ;q,t\, ]
\Big(\prod_{i=1}^k (1-\, t^{k-i}q^{\la_i})\Big)P_\la [Y_k;q,t\, ]
\eqno 1.5
$$
Now it follows from theorem (4.17) of [12] (p.~325) that if $l(\la)=k$ then
$$
P_\la [Y_k;q,t\, ]=y_1y_2\cdots y_k\, P_\mu [Y_k;q,t\, ] 
$$
with
$$
\mu=(\la_1-1,\la_2-1,\ldots ,\la_k-1)\ess .
$$
Thus we may rewrite 1.5 in the form
$$
M^{(Y_k)}(-1)\OM\big[X_nY_k{\textstyle{1-t\over 1-q}}\big]= y_1y_2\cdots y_k \sum_{l(\mu)\leq k} Q_{\mu+1^k}[X_n ;q,t\, ]
\Big(\prod_{i=1}^k (1-\, t^{k-i}q^{\mu_i+1})\Big)P_\mu [Y_k;q,t\, ]
\eqno 1.6
$$
Now from I.5 and I.6 we derive  that
$$
Q_{\mu+1^k}[X_n ;q,t\, ]\Big(\prod_{i=1}^k (1-\, t^{k-i}q^{\mu_i+1})\Big)\ses {J_{\mu+1^k}[X_n ;q,t\, ]\over h_\mu'(q,t)
}\,\,  .
$$
Thus the desired identity in 1.4 immediately follows upon  using this in 1.6.
\sa

\noindent{\bol Proposition 1.3}
$$
{M^{(Y_k)}(u)\OM\left[  X_n Y_k{\textstyle{1-t\over 1-q}}\right] \over \OM\left[  X_n Y_k{\textstyle{1-t\over 1-q}}\right]}\ses
\sum_{r=0}^k  u^r\sum_{I\con [1,k] \atop |I|=r}   
A_I[Y_k;t]\,\, \OM\left[ X_n  Y_I  (t-1) \right]\, ,
\eqno 1.7
$$
{\ita in particular this ratio is independent of $q$.}
\sas

\noindent{\bol Proof}
$$
\eqalign{
M^{(Y_k)}(u)\OM\left[  X_n Y_k{\textstyle{1-t\over 1-q}}\right]
&\ses  
\sum_{r=0}^k  u^r\sum_{I\con [1,k]\atop |I|=r}   A_I[Y_k;t] \OM\left[  X_n
\big(Y_k+(q-1)Y_I \big){\textstyle{1-t\over 1-q}}\right]\cr &\ses \OM\left[ 
X_n Y_k{\textstyle{1-t\over 1-q}}\right]
\sum_{r=0}^k  u^r\sum_{I\con [1,k]\atop |I|=r}   A_I[Y_k;t]\, \OM\left[  X_n 
Y_I  (t-1) \right]
\cr
}
$$
and this proves 1.7.
\sa

\noindent{\bol Theorem 1.1}

{\ita If $B_k(X;t)$ is any operator with the property that
$$
B_k(X;t) J_\mu[X;0,t]\ses J_{\mu+1^k}[X;0,t]
\ess\ess\ess\ess \hbox{for all $l(\mu)\leq k$}
\eqno 1.8
$$ 
then }
$$
\htq{B_k}(X;t)J_\mu[X;q,t]\ses J_{\mu+1^k}[X;q,t]
\ess\ess\ess\ess \hbox{for all $l(\mu)\leq k$}
\eqno 1.9
$$ 
{\ \bol Proof}

Suppose we show that
$$
y_1y_2\cdots y_k \htq{ B_k}(X_n;t)^{X_n}\OM\left[  X_n Y_k{\textstyle{1-t\over
1-q}}\right]\ses M^{(Y_k)}(-1)\OM\left[  X_n Y_k{\textstyle{1-t\over
1-q}}\right]\ess ,
\eqno 1.10
$$
then by combining this with Proposition 1.2 we get that
$$
\eqalign{
y_1y_2\cdots y_k\, \sum_{l(\mu)\leq k} &  \htq{B_k}(X_n;t)^{X_n}J_\mu[X_n;q,t]
{P_\mu[Y_k;q,t]\over h_\mu'(q,t)} 
\cr
&\ess\ess\ess\ess\ses
y_1y_2\cdots y_k\,  
\sum_{l(\mu)\leq k }J_{\mu+1^k}[X_n;q,t]\,{ P_\mu[Y_k;q,t]\over h_\mu'(q,t)}
\cr}
$$
and 1.9 then follows by equating coefficients of $P_\mu[Y_k;q,t]$.
To show 1.10 we need only verify that
$$
y_1y_2\cdots y_k {\htq{ B_k}(X_n;t)^{X_n}\OM\left[  X_n
Y_k{\textstyle{1-t\over 1-q}}\right]
\over \OM\left[  X_n Y_k{\textstyle{1-t\over 1-q}}\right] }\ses
{ M^{(Y_k)}(-1)\OM\left[  X_n Y_k{\textstyle{1-t\over 1-q}}\right] \over \OM\left[  X_n Y_k{\textstyle{1-t\over 1-q}}\right] } \ess ,
\eqno 1.11
$$
and since we have shown that both sides of this equation are independent of $q$, we need only verify
this equality at $q=0$.  However, the hypothesis in 1.8 yields that
$$
\eqalign {
y_1y_2\cdots y_k{\htq{ B_k}(X_n;t)^{X_n}\OM\left[  X_n Y_k{\textstyle{1-t\over
1-q}}\right]
\over \OM\left[  X_n Y_k{\textstyle{1-t\over 1-q}}\right] }\bigg|_{q=0}
&\cr
={y_1y_2\cdots y_k\over  \OM\left[  X_n Y_k{\textstyle{(1-t) }}\right]} 
&\ess 
\sum_{l(\mu)\leq k} B_k(X_n;t)^{X_n}J_\mu[X_n;0,t] \,\,
{ P_\mu[Y_k;0,t]\over h_\mu'(0,t)}\cr
&=
{y_1y_2\cdots y_k\over  \OM\left[  X_n Y_k{\textstyle{(1-t) }}\right]}
\sum_{l(\mu)\leq k} J_{\mu+1^k}[X_n;0,t] \,\,
{ P_\mu[Y_k;0,t]\over h_\mu'(0,t)}\cr
}
$$
and again by Proposition 1.3 we see that this is precisely
$$
{ M^{(Y_k)}(-1)\OM\left[  X_n Y_k{\textstyle{1-t\over 1-q}}\right] \over \OM\left[  X_n Y_k{\textstyle{1-t\over 1-q}}\right] }
\bigg|_{q=0}\,\,  .
$$
This completes the proof of 1.9.
\sa

Before we can proceed with the proof of Theorem I.1 we need one more auxiliary result. 
To begin with it will be convenient to consider the substitution $X\RA X/(1-t)$ as a linear
operator on symmetric functions. More precisely, for any symmetric polynomial $P$ and any alphabet $X$
we set
$$
F^tP[X]\ses P\big[{\textstyle{X(1-t)}}\big]
$$ 
Now a somewhat surprising development is that 
the operation $V\RA \htq{V}$  defined in I.11  commutes with conjugation by
$F^t$. In fact, we may state
\sas

\noindent{\bol Proposition  1.4}

{\ita 
For any linear operator $V$ acting on $\Lambda$ and any polynomial $P\in \Lambda$
we have}
$$
F^t\, \htq{V} {F^t}^{-1}P[X]\ses \htq{(F^t\,  {V} {F^t}^{-1})} P[X]
\eqno 1.12
$$
{\bol Proof}

It is sufficient to prove 1.12 for the Schur basis.  Note that
for any partition $\la$, 
 the addition formula for Schur functions gives
$$ 
\eqalign{
\htq{(F^t{V}{F^t}^{-1})} S_\la[X] &\ses (F^t\,  {V} {F^t}^{-1})^Y
S_\la[qX+(1-q)Y]\Big|_{Y=X}
\cr
&\ses \sum_{\mu\con \la} S_{\la/\mu}[qX]\, (F^t\,  {V} {F^t}^{-1})^Y S_\mu[(1-q)Y]\Big|_{Y=X}
\cr
&\ses\sum_{\mu\con \la} S_{\la/\mu}[qX]\, (F^t\,  {V} {F^t}^{-1})^X S_\mu[(1-q)X]\ess .
\cr
}
\eqno 1.13
$$
In the same vein we see that the left hand side of 1.12, for $P=S_\la$,   gives
$$
\eqalign{
F^t\,\htq{V} {F^t}^{-1}S_\la[X] &\ses F^t\,  \htq{V}
S_\la[{\textstyle{X\over 1-t}}]
\cr
&\ses  F^t \, V ^Y S_\la \left[{\textstyle{qX+(1-q)Y\over 1-t}} \right]|_{Y=X}
\cr
&\ses  F^t\, \sum_{\mu\con \la} S_{\la/\mu}\big[{\textstyle{qX \over 1-t}}\big]\,  V ^Y S_\mu\big[{\textstyle{  {(1-q) Y\over 1-t}
}}\big]\Big|_{Y=X}
\cr
&\ses  \sum_{\mu\con \la} S_{\la/\mu}\big[{\textstyle{qX  }}\big]\,  F^t\, V ^X S_\mu\big[{\textstyle{  {(1-q) X\over 1-t}
}}\big] 
\cr
}  
$$ 
and it is easily seen that this is another way to write the last expression in 1.13.
\sa

We are now in a position to give our
\sas

\noindent
{\bol Proof of Theorem I.1}

By assumption we have 
$$
\OH_k H_\mu[X;t\, ]\ses H_{\mu+1^k} [X;t\, ]
\ess\ess\ess\ess \hbox{for all $\mu$ of length $\leq k$} 
\eqno 1.14
$$
Thus from I.10 we derive that
$$
\OH_k{F^t}^{-1}Q_\mu[X;t\, ]\ses  {F^t}^{-1}Q_{\mu+1^k} [X;t\, ]\, .
$$
Now this, using I.9, may be rewritten as
$$
F^t \OH_k{F^t}^{-1}J_\mu[X;0,t\, ]\ses J_{\mu+1^k} [X;0,t\, ]\ess\ess\ess\ess
\hbox{for all $\mu$ of length $\leq k$} \, .
$$
In other words the operator
$$
B_k=B_k(X;t)\ses F^t\, \OH_k{F^t}^{-1}
$$
satisfies the hypothesis of Theorem 1.1. It then  follows that
$$
{\htq{({F^t \OH_k{F^t}^{-1}})}}J_\mu[X;q,t\, ]\ses J_{\mu+1^k} [X;q,t\, ]
\ess\ess\ess\ess \hbox{for all $\mu$ of length $\leq k$} \,
$$
But Proposition 1.4 yields 
$$
{F^t {\htq{\OH_k} }}{F^t}^{-1}J_\mu[X;q,t\, ]\ses J_{\mu+1^k} [X;q,t\, ]
\ess\ess\ess\ess \hbox{for all $\mu$ of length $\leq k$} \,
$$
and I.7 shows that this is just another way of writing I.13,
completing the proof of  Theorem I.1.
\sa

\vbox{
\noindent{\bol Remark 1.1}

We should note that any  symmetric function 
operator $V(q,t)=V(x;q,t)$  of the form
$$
V(q,t)\ses \sum_{I\con [1,k]}a_I(x,t)\ssp T_I^q\ess  .
\eqno 1.15
$$
(in particular the Macdonald operator)
satisfies the identity
$$
 V(q,t)\ses  \htq{V (0,t)}  
\eqno 1.16
$$}

\noindent
In fact, since for any $P\in\Lambda$ we have
$$
V(q,t)P[X]\ses \sum_{I\con [1,k]}a_I(x,t)\ssp P[X+(q-1)X_I]\, .
$$
then 
$$
V(0,t)\, P[X] \ses\sum_{I\con
[1,k]}a_I(x,t)\ssp P[X-X_I]\, .
$$
Thus  
$$
\eqalign{
\htq{V(0,t)}  P[X]&\ses V(0,t)^Y P[qX+(1-q)Y ]\, \Big|_{Y =X}
\cr
&\ses \sum_{I\con [1,k]}a_I(y,t)\ssp P[qX+(1-q)(Y -Y_I)]\, \Big|_{Y =X}
\cr
&\ses \sum_{I\con [1,k]}a_I(x,t)\ssp P[qX+(1-q)(X-X_I)] 
\cr
&\ses \sum_{I\con [1,k]}a_I(x,t)\ssp P[X+(q-1)X_I] 
\cr
}
$$
which is 1.16. 
\sa 

We then see that Propositions 1.1 and 1.3 are both particular
cases of the following general fact:
\sas

\noindent{\bol Proposition 1.5}  
\sas

{\ita If $V(q,t)=V(x;q,t)$ is an operator on $\Lambda$ with the property
$$
V( q,t)\ses \htq{V( 0,t)} 
$$
then
$$
{V( q,t)^X\OM\big[XY{\textstyle{1-t\over 1-q}}\big]\over \OM\big[XY{\textstyle{1-t\over 1-q}}\big]}\ses 
 { V( 0,t)^X\OM\big[XY{ (1-t)}\big]\over \OM\big[XY{ (1-t)}\big] }
$$
in particular this ratio is independent of $q$.
}

The proof follows exactly the same steps we used to prove Proposition 1.1.

\sap
\vbox{
\noindent{\bol 2. Applications}

It is shown in [12] ((4.14) p.~324) that
$$
P_\mu[X;t,t]\ses S_\mu[X]\ess .
\eqno 2.1
$$ 
thus I.5 gives 
$$
J_\mu[X;t,t]\ses \prod_{s\in \mu}(1-t^{h_\mu(s)}) S_\mu[X]\ess ,
\eqno 2.2
$$
}\noindent
where $h_\mu(s)=a_\mu(s)+l_\mu(s)+1$ denotes the hook-length corresponding
to the cell $s$ in $\mu$. This given, Theorem 1.1 has the following immediate
corollary.
\sas

\noindent{\bol Theorem 2.1}

{\ita If  $B_k$ is any operator on $\Lambda$ with the property 
$$
B_k(X;t)J_\mu[X;0,t]\ses J_{\mu+1^k}[X;0,t]
\ess\ess\ess\ess \hbox{for all $l(\mu)\leq k$}
\eqno 2.3
$$ 
then  the operator $\htt{B_k}$ defined by setting for $P\in \Lambda$
$$
\htt{ B_k}(X;t)P[X]\ses B_k(Y;t)^YP[tX+(1-t)Y]\big|_{Y=X}
\eqno 2.4
$$
has the property}
$$
\htt{B_k}(X;t)S_\mu[X]\ses \Big(\prod_{i=1}^k\big
(1-t^{k+1-i+\mu_i}\big)\Big) \ssp S_{\mu+1^k}[X]
\ess\ess\ess\ess \hbox{for all $l(\la)\leq k$}\ess .
\eqno 2.5
$$ 
{\ \bol Proof}

This follows by setting $q=t$ in 1.9, using formula 2.2 and canceling the common factor.
\sa

This result has the following  converse
\sas

\noindent{\bol Theorem 2.2}

{\ita If $B_k(X;t)$ is any operator with the property in 2.5,
then  the operator $\htq{ B_k}(X,t)$ defined by setting for $P\in
\Lambda$,
$$
\htq{B_k}(X;t)P[X]\ses B_k(Y;t)^YP[qX+(1-q)Y]\big|_{Y=X}
\eqno 2.6
$$
has the property}
$$
\htq{B_k}(X;t)J_\mu[X;q,t]\ses J_{\mu+1^k}[X;q,t]
\ess\ess\ess\ess \hbox{for all $l(\mu)\leq k$}
\eqno 2.7
$$ 
 
\noindent{\bol Proof}

The argument follows almost verbatim what we did to prove Theorem 1.1 except that
at one point we must set $q=t$ rather than $q=0$. 
\sas

Of course we may produce a polynomiality proof based on this theorem, however
all our attempts yielded a more complicated argument than that based on Theorem I.1. 
\sa

For representation theoretical reasons Garsia-Haiman where led to consider the polynomials
$$
\TH_\mu[X;q,t]\ses \sum_{\la} S_\la[X] \TK_{\la\mu}(q,t)
\eqno 2.8
$$
where
$$
\TK_{\la\mu}(q,t)\ses t^{n(\mu)}K_{\la\mu}(q,1/t)
\eqno 2.9
$$
with
$$
n(\mu)=\sum_i(i-1)\mu_i
\eqno 2.10
$$
Note that setting $q=0$ in 2.8 gives
$$
\TH_\mu[X;0,t]\ses \TH_\mu[X;t]=\sum_{\la} S_\la[X]\TK_{\la\mu}(t)
$$
where 
$$
\TK_{\la\mu}(t)\ses t^{n(\mu)} K_{\la\mu}(1/t)
$$
is the so-called {\ita cocharge} Kostka-Foulkes polynomial. Now it develops
that a result analogous to Theorem I.1 holds also for the basis
$\big\{\TH_\mu[X;q,t]\big\}_\mu$.
\sas
\def \TTH {{\cal H}}
\noindent{\bol Theorem 2.3}

{\ita
If $\TTH_k=\TTH_k(X,t)$ is any linear operator on $\Lambda$ with the property
that
$$
\TTH_k \TH_\mu[X;t\, ]\ses \TH_{\mu+1^k} [X;t\, ]
\ess\ess\ess\ess \hbox{for all $\mu$ of length $\leq k$} 
\eqno 2.11 
$$
then the operator  $\htq{\TTH_k}$ defined by I.11 has the property
$$
\htq{\TTH_k} \TH_\mu[X;q,t\, ]\ses \TH_{\mu+1^k} [X;q,t\, ]
\ess\ess\ess\ess \hbox{for all $\mu$ of length $\leq k$}\ess . 
\eqno 2.12 
$$
In particular, the modified Macdonald polynomial $\TH_\mu[X;q,t\,]$ may be
obtained from the ``Rodriguez'' formula:
$$
\TH_\mu[X;q,t\,]\ses \htq{\TTH_{\mu_1'}}\htq{\TTH_{\mu_2'}}\cdots
\htq{\TTH_{\mu_h'}}\, {\bf 1}
\eqno  2.13
$$
where $\mu'=(\mu_1',\mu_2',\ldots ,\mu_h')$ denotes the conjugate of $\mu$.
}

\noindent{\bol Proof}

Note that setting $t=1/t$ in 2.11 gives 
$$
t^{-n(\mu)} \TTH_k(X;1/t)  H_\mu[X;t\, ]\ses t^{-n(\mu+1^k)}  H_{\mu+1^k} [X;t\,
]\ssp .
\ess\ess\ess\ess \hbox{for all $\mu$ of length $\leq k$}  
$$
and this (using 2.10)  may be rewritten as 
$$
t^{{k \choose 2}}\TTH_k(X;1/t)  H_\mu[X;t\, ]\ses    H_{\mu+1^k} [X; t\, ]
\ess\ess\ess\ess \hbox{for all $\mu$ of length $\leq k$} \ssp . 
$$
Thus we may apply Theorem I.1 and derive that
$$
t^{{k \choose 2}}\TTH_k(Y;1/t)^Y \,  H_\mu[qX+(1-q)Y;q,t\, ]\big|_{Y=X}= 
 H_{\mu+1^k}[ X ;q,t\, ]
\ess\ess\ess\ess \hbox{for all $\mu$ of length $\leq k$}\ssp .
$$
Setting $t=1/t$ and using 2.9 transforms this back to 2.12.
\sa

Some interesting developments follow by combining the present
identities with one of the 
simplest of the ``Rodriguez'' operators introduced by Lapointe-Vinet in [10].
To see how this comes about we need to recall 
this beautiful result.
\sas

\noindent{\bol Theorem 2.4}

{\ita The operator 
$$
LV_k(X_n;q,t) \ses {1\over ({1\over q};{1\over t})_{n-k}}M^{(X_n)}
\big(-{\textstyle {1\over qt^{n-k-1}}}\big)\underline {e_k[X_n]}
\eqno 2.14
$$
where $\underline {e_k[X_n]}$ denotes multiplication by $ e_k[X_n] $, has
the property
$$
LV_k(X_n;q,t)J_\mu[X_n;q,t\, ]\ses J_{\mu+1^k}[X_n;q,t\, ]
\ess\ess\ess\ess\ess\hbox{for all $\ess l(\mu)\leq k$}
$$
}
\noindent{\bol Proof}

The argument is so elementary that it might as well be reproduced here.
It is shown by Macdonald in [12] (p.~340 (6.24) (iv)) that
$$
e_k[X_n]P_\mu[X_n;q,t\,] \ses P_{\mu+1^k}[X_n;q,t\, ]\sps 
\sum_{\la/ \mu\in V_k
\atop \la \neq \mu+1^k} P_\la[X_n;q,t\,] \Psi_{\la\mu}(q,t)\ess
$$ 
where $\la/ \mu\in V_k$ indicates that the sum is over partitions $\la$ such
that
$\la/\mu$ is a vertical $k$-strip. 
Now applying the Macdonald operator $M^{(X_n)}(u)$ to both sides and using
(4.15) p.~324 of [12]  we get for $l(\mu)\leq k$
$$
\eqalign{
M^{(X_n)}(u)
e_k[X_n]P_\mu[X_n;q,t\,] 
&=\prod_{i=1}^k\big(1+u\, t^{n-i}q^{\mu_i+1}\big) \prod_{i=k+1}^n
\big(1+u\, t^{n-i}\big)  \ess
P_{\mu+1^k}[X_n;q,t\,] \cr &\ess\ess\ess\ess\ess\ess\ess\ess\ess\ess 
\sps \sum_{\la/ \mu\in V_k
\atop \la \neq \mu+1^k}\ess \prod_{i=1}^n
\big(1+u\,t^{n-i}q^{\la_i}\big) P_\la[X_n;q,t\,]
\Psi_{\la\mu}(q,t)\ssp .
\cr
}
\eqno 2.15
$$
Now $1^k$  is the shortest vertical $k$-strip  that may be added to $\mu$, and  any other will spill a 
cell at height $k+1$. Thus each term of the sum in 2.15 will contain the factor
$$
\big(1+u\,t^{n-k-1}q\big)
$$  
which vanishes when we set $u=-1/q\, t^{n-k-1}$. Thus 2.15 gives
$$
{1\over \big({\textstyle  {1\over q};{1\over t}} \big)_{n-k}}\ssp M^{(X_n)}\big
({-\textstyle {1\over qt^{n-k-1}}}\big) e_k[X_n]P_\mu[X_n;q,t\,] 
\ses  \prod_{i=1}^k\big(1- \, t^{k+1-i}q^{\mu_i }\big)P_{\mu+1^k}[X_n;q,t\, ]
$$
and this is easily converted to 2.14 by means of I.5.
\sas

To state and prove our next result we need some notation. To begin, if $\sig=(\sig_1,\sig_2,\ldots ,\sig_n)$ is 
a permutation in the
symmetric group ${\cal S}_n$ and $a=(a_1,a_2,\ldots, a_n)$ is a given vector,  
we set
$$
\sig a\ses (a_{\sig_1},a_{\sig_1},\ldots ,a_{\sig_n})\ssp ,
$$ 
we also let $Supp(a)$ denote the ``{\ita support}'' of $a$, that is the set of
elements
$$
Supp(a)\ses \{i\, :\, a_i\neq 0\, \}\ssp .
$$
Next, for any two vectors  $a=(a_1,a_2,\ldots, a_n)$ and $b=(b_1,b_2,\ldots, b_n)$ we shall write $a\approx b$ if and only if
the components of $b$ are a rearrangement of the components of $a$. More precisely, we set $a\approx b$ if
and only if for some
$\sig\in {\cal S}_n$ we have
$$
b\ses  \sig \, a
$$ 
With this notation the monomial symmetric function $m_\la$ may be represented by the sum
$$
m_\la[X_n]\ses \sum_{p\,\,  \approx \la}\, x^p\ssp .
\eqno 2.16
$$
Finally, we shall generically denote by $\eee=(\eee_1,\eee_2,\ldots ,\eee_n)$ the indicator vector  of
a subset of $\{1,2,\ldots ,n\}$. Thus setting $|\eee|=k$ will mean that $\eee$ represents a subset of cardinality $k$. 
In particular, the elementary symmetric function $e_k$  may be represented by the sum
$$
e_k[X_n]\ses \sum_{|\eee|=k}\ssp x^\eee
\eqno 2.17
$$
\sas

It develops that, not withstanding the presence of terms $1/q$ in 2.14, we can
evaluate the limit of the operator $LV_k(X_n;q,t)$ as $q\RA 0$. What follows
is the following surprising corollary of the Lapointe-Vinet result. 
\sas

\noindent{\bol Theorem 2.5}

{\ita Let $TLV_k(t)$ be the operator defined by setting for the monomial basis
$$
TLV_k(t)\, m_\la\ses 0\ess\ess\ess \hbox{when  $l(\la)>k$}
\eqno 2.18
$$
and
$$
TLV_k(t)\, m_\la\ses \sum_{|\eee|=k}\ssp \sum_{ {p\approx \la  }  }\,c_{\eee,p}(t)\,  S_{p+\eee} 
\ess\ess\ess \hbox{when  $l(\la)\leq k$}
\eqno 2.19
$$
where $S_{p+\eee}$ denotes the corresponding signed Schur function and 
$$
c_{\eee,p}(t)\ses 
\cases {
(-1)^{n-k}t^{n-k \choose 2} \prod_{\eee_i=0}\big(-t^{k+1-i}\big)\prod_{\eee_i+ p_i=1}\big(1-t^{k+1-i}\big) & if $ Supp(p)\con
Supp(\eee)$
\cr\cr
0 & otherwise
\cr
}
\eqno 2.20
$$
then we also have } 
$$
TLV_k Q_\mu[X;t\, ]\ses 
\cases{
Q_{\mu+1^k} [X;t\, ]& if  $l(\mu)\leq k$\cr\cr 
0 & otherwise
}
\eqno 2.21
$$
{\bol Proof}

The original definition, of the Macdonald operator, given in (3.2) p.~315 of
[12] may be written in the form
$$
M^{(X_n)}(u)\ses {1\over \Delta_n(x)}\sum_{\sig\in {\cal S}_n}sign(\sig)
x^{\sig\delta }
\prod_{i=1}^n \Big(1+ut^{n-\sig_i}T_i^q  \Big)
\eqno 2.22
$$
where $\Delta_n(x)$ denotes the Vandermonde determinant in the variables $\xon$, 
$T_i^q$ denotes the operator that replaces $x_i$ by $qx_i$ and for convenience we have
set 
$$ 
\delta\ses 
(n-1,n-2,\ldots ,1,0)\ssp .
$$ 
This given, 
taking account of 2.16 and 2.17 we may write the action of the Lapointe-Vinet operator $LV_k$ 
on the monomial basis in the form
$$
LV_k\, m_\la[X_n]= {1\over ({1\over q};{1\over t})_{n-k}}
\sum_{|\eee|=k}\, \sum_{p\,\, \approx \la }
\ssp 
{1\over \Delta_n(x)}\sum_{\sig\in {\cal S}_n}sign(\sig)
x^{\sig\delta }
\prod_{i=1}^n\Big(1-t^{k+1-\sig_i}q^{ p_{\sig_i}+\eee_{\sig_i}-1 }\Big)x^{\sig (p+\eee)}\, .
\eqno 2.23
$$
Since we may set
$$
{1\over \Delta_n(x)}\sum_{\sig\in {\cal S}_n}sign(\sig)
x^{\sig(\delta+p+\eee) }\ses S_{p+\eee}[X_n]
$$
we see that 2.23 reduces to
$$
LV_k(q,t)\, m_\la[X_n] \ses \sum_{|\eee|=k}\, 
\sum_{p\,\, \approx \la } c_{\eee,p}(q,t)\,\,  S_{p+\eee}[X_n]
\eqno 2.24
$$
with
$$
c_{\eee,p}(q,t)\ses {1\over ({1\over q};{1\over t})_{n-k}}
\prod_{i=1}^n\big(1-t^{k+1- i}q^{ p_{ i}+\eee_i-1 }\big)
\ssp .
\eqno 2.25
$$
Our next step is to evaluate  $c_{\eee,p}(q,t)$ at $q=0$. To this end
it is convenient to rewrite 2.25 in the form
$$
\eqalign{
c_{\eee,p}(q,t)&= {(-1)^{n-k}t^{n-k \choose 2}\over ({  q};{  t})_{n-k}}
{q^{n-k} \over q^{\#\{i\, :\, \eee_i+p_i=0\}} }
\prod_{\eee_i+p_i>1}\big(1-t^{k+1- i}q^{ p_{ i}+\eee_i-1 }\big)\ess  \cr
&\ess\ess\ess\ess\ess\ess\ess\ess\ess\ess\ess\ess\ess\ess\ess\ess\ess\ess\ess\ess\ess\ess\ess\ess\ess
\times
 \prod_{\eee_i+p_i=1} \big(1-t^{k+1- i}\big) 
 \prod_{\eee_i+p_i=0}\big(q-t^{k+1- i} \big) 
\cr
}
\eqno 2.26
$$
Now note that since $|\eee|=k$ we necessarily have that
$$
\#\{i\, :\, \eee_i+p_i=0\}\leq n-k
$$
with equality only if 
$$
Supp(p)\con Supp(\eee)\ssp .
\eqno 2.27
$$
Thus it follows that 
$$
{q^{n-k} \over q^{\#\{i\, :\, \eee_i+p_i=0\}} }\, \Big|_{q=0}\ses 0
$$
when 2.27 fails to hold.
Since we have the obvious evaluations
$$
\eqalign{
  ({  q};{  t})_{n-k}\,\Big|_{q=0}  &\ssp =\,  1
\cr
 \prod_{\eee_i+p_i>1}\big(1-t^{k+1- i}q^{ p_{ i}+\eee_i-1 }\big)\Big|_{q=0} & \ssp  =\,  1
\cr
 \ess  \prod_{\eee_i+p_i=0}\big(q-t^{k+1- i} \big) \,\Big|_{q=0} & \ssp  =\,  \prod_{\eee_i+p_i=0}\big( -t^{k+1- i} \big)
\cr}
$$
we see that from 2.26 we derive that
$$
 \, c_{\eee,p}(0,t)\, \ses  c_{\eee,p}(t)\, ,
$$
with $ c_{\eee,p}(t)$ precisely as defined in 2.20. Thus  from 2.24 we get
$$
LV_k(0,t)\, m_\la[X_n]\ses \sum_{|\eee|=k}\, \sum_{p\,\, \approx \la }
c_{\eee,p}(0,t)\,\,  S_{p+\eee}[X_n]
\ses TLV_k(t)\, m_\la[X_n]\ess .
$$
To show 2.21 note that from I.8 we get
$$
LV_k(q,t) J_\mu[X_n;q,t\, ]\,  \ses\sum_{\la}
\Big(LV_k(q,t)S_\la[X_n(1-t)]\Big)\,  K_{\la\mu}(q,t)\ssp .
$$
Thus the polynomiality of the $ K_{\la\mu}(q,t)$ assure that we can safely set $q=0$ here and
obtain that
$$
\eqalign{
LV_k(q,t)J_\mu[X_n;q,t\, ]\, \Big|_{q=0}&\ses\sum_{\la}
\Big(LV_k(0,t)S_\la[X_n(1-t)]\Big)\,  K_{\la\mu}(0,t)
\cr
&\ses\sum_{\la} \Big(TLV_k(t)S_\la[X_n(1-t)]\Big)\,  K_{\la\mu}(t)
\cr
&\ses TLV_k(t) Q_\mu[X_n;t\, ]
\cr
}
$$
Thus when $l(\mu)\leq k$ the Lapointe-Vinet result (Theorem 2.5) yields
$$
TLV_k(t) Q_\mu[X;t\, ]\ses J_{\mu+1^k}[X;0,t\, ]\ses Q_{\mu+1^k}[X;t\,
]\ssp .
$$
This proves the first alternative in 2.21. To show the second we note that (2.6) of
[12] p.~209 implies that we have an expansion of the form
$$
J_\mu[X;t]=\sum_{\la\leq \mu} m_\la[X] \, \xi_{\la\mu}(t)\ssp .
\eqno 2.28
$$ 
Thus when $l(\mu)>k$ we shall have $l(\la)  >k$  for all the summands in 
2.28 and 2.18 then gives that
$$
TLV_k(t)\,J_\mu[X_n;t] 
\ses 0\ssp .
$$ 
This completes our proof.
\sa
 
We should note that the Lapointe-Vinet result has one further
curious consequence. 
\sa

\noindent{\bol Theorem 2.6}

{\ita
Let $W_k=W_k(X_n;t)$ be the operator defined by setting for every $P\in \Lambda$
$$  
W_k P[X]\ses {1\over ({1\over t};{1\over t})_{n-k}}\sum_{r=0}^n \Big(
{-1\over t^{n-k}}\Big)^r
\sum_{I\con [1,n]\atop |I|=r}B_I(x;t)\, P[X-X_I]
\eqno 2.29
$$ 
with
$$
B_I(x;t)\ses {1\over \Delta_n(x)} T_I^t\, \Delta_n(x)e_k[X_n]\ssp .
\eqno 2.30
$$
Then}
$$
\htq{W_k } \, J_\mu[X_n;q,t\, ]\ses  J_{\mu+1^k}[X_n;q,t\, ]
\ess\ess\ess\ess\ess\ess
\hbox{for all $\mu$ of length $\leq k$.}
\eqno 2.31
$$
{\bol Proof}

Using formula (3.5) of [12] p.~316 we may write the Lapointe-Vinet result 
\hbox{(for $l(\mu)\leq k$)} in the form
$$  
J_{\mu+1^k}[X_n;q,t\, ]= 
{1\over ({1\over q};{1\over t})_{n-k}}\sum_{r=0}^n
\Big( {-1\over q\,  t^{n-k-1}}\Big)^r
\hskip -.1in\sum_{I\con [1,n]\atop |I|=r}\hskip -.08in { T_I^t\,
\Delta_n(x)\over
\Delta_n(x)}\, T_I^q e_k[X_n] J_{\mu}[X_n;q,t\, ]
$$
where, for each $i\in I$, $T_I^t$ replaces $x_i$ by $tx_i$.
This given, we may set $q=t$ and obtain that
$$  
\eqalign{
J_{\mu+1^k}[X_n;t,t\, ]&=  {1\over ({1\over t};{1\over t})_{n-k}}\sum_{r=0}^n 
\Big( {-1 \over t^{n-k}}\Big)^r
 \sum_{I\con [1,n]\atop |I|=r}  B_I(x;t)\,
\,  J_{\mu}[X_n+(t-1)X_I;t,t\, ]
\cr
&= \htt{W_k} \, J_\mu[X_n;t,t\, ]
\ssp .
}
\eqno 2.32
$$ 
Since
$$
J_{\mu}[X_n;t,t\, ]\ses h_\mu(t) S_\mu[X_n]\ssp ,
$$
we derive from 2.32 that the operator $W_k $ satisfies the hypothesis of
Theorem 2.2, thus 2.31 is an immediate consequence of  Theorem 2.2.
\sas

\vbox{
We terminate with one final application of Theorem I.1:
\sas

\noindent{\bol Theorem 2.7}

{\ita Let $\OH_k(X;t)$ be any operator with the property that
$$
\OH_k(X;t)H_\mu[X;t]\ses H_{\mu+1^k}[X;t]\ess\ess\ess\ess\hbox{for all $\mu$ of
length $\ssp \leq k$}.
\eqno 2.33
$$
Then the operator
$$
\htt{\OH_k}(X;q)\ses \htq{\OH_k}(X;t)\, \Big|_{t\,\,  \rightleftharpoons \,\, q} 
$$ 
has the property
$$
\omega \htt{\OH_k}(X;q) \omega H_\mu[X;q,t\, ]\ses H_{(k,\mu)}[X;q,t\, ]
\ess\ess\ess\ess\hbox{for all $\mu$ with  $\ssp \mu_1 \leq k$}.
\eqno 2.34
$$
}

\noindent{\bol Proof}

It follows from the Macdonald duality formula [12] ((5.1) p.~327) that the
polynomial $H_\mu[X;q,t\, ]$ satisfies the identity
$$
H_{\mu'} [X;q,t\, ]  \ses \omega \, H_\mu[X;t,q\, ]
\eqno 2.35
$$
Now, assuming 2.33, from Theorem I.1 it follows that
$$
 \htq{\OH_k}(X;t) \,H_{\mu'}[X;q,t\, ]\ses H_{\mu'+1^k}[X;q,t\, ]
\ess\ess\ess\ess\hbox{for all $\mu'$ of length $\ssp \leq k$}. 
$$
Interchanging $q$ and $t$ we get
$$
 \htt{\OH_k}(X;q)\,H_{\mu'}[X;t,q\, ]\ses H_{\mu'+1^k}[X;t,q\, ]
\ess\ess\ess\ess\hbox{for all $\mu'$ of length $\ssp \leq k$}. 
$$
and two uses of 2.35 then give  
$$
 \htt{\OH_k}(X;q)\,\omega H_{\mu }[X;t,q\, ]\ses \omega  H_{(k,\mu)}[X;q,t\, ]
\ess\ess\ess\ess\hbox{for all $\mu $ with  $\ssp \mu_1 \leq k$}, 
$$
which simply another way of writing formula 2.34.

\sas

The simplicity of our operator $\omega \htt{\OH_k}(X;q)\,\omega$ and our proof
of 2.34 should be contrasted with the complexity of the developments in  the
recent 
 Kajihara-Noumi paper [5]. 
}

\vfill\supereject
\centerline{\bol REFERENCES} 
\sas
\item{[1]}
A.~Garsia and M.~Haiman,
  {\ita Some natural bigraded ${S}_n$-modules and $q,t$-{K}ostka
  coefficients}, Electron. J. Combin. {\bol 3} (1996), no.~2, Research Paper 24,
  approx.\ 60 pp.\ (electronic),
  The Foata
  Festschrift,
\sas

\item{[2]}
A. Garsia, M. Haiman and Tesler, {\ita Explicit Plethystic Formulas 
for Macdonald $q,t$-Kostka Coefficients }, The Andrews Festschrift,
S\'eminaire Lotharingien de Combinatoire {\bol 42},  B42m.
Website http://www.emis.de/journals/SLC/.
\sas 

\item{[3]}
A.~M.\ Garsia and J.~Remmel, {\ita {P}lethystic {F}ormulas and positivity for
  $q,t$-{K}ostka {C}oefficients}, {M}athematical {E}ssays in {H}onor of
  {G}ian-{C}arlo {R}ota (B.~E.\ Sagan and R.~Stanley, eds.), Progress in
  Mathematics, vol.~161, 1998.
\sas 

\item{[4]}
A.~M.\ Garsia and G.~Tesler, {\ita Plethystic formulas for {M}acdonald
  $q,t$-{K}ostka coefficients}, Adv.\ Math.\ {\bol 123} (1996), no.~2, 144--222.
\sas

\item{[5]} Y. Kajihara and M. Noumi, {\ita Raising Operators of Row Type for
Macdonald Polynomials}, Compositio Mathematica {\bf 120} (2000), 119-136.
\sas

\item{[6]} A. ~Kirillov and M. Noumi, {\ita q-difference raising operators for Macdonald polynomials
and the integrality of transition coefficients}, q-alg/9605005. 4$\,$ May 96.
\sas

\item{[7]}
F.~Knop, {\ita Integrality of two variable {K}ostka functions}, J.\ Reine
  Angew.\ Math. {\bol 482} (1997), 177--189.
\sas

\item{[8]}
F.~Knop, {\ita Symmetric and non-symmetric quantum {C}apelli polynomials},
  Comment. Math. Helv. {\bol 72} (1997), no.~1, 84--100.

\sas
\item{[9]}
L.~Lapointe and L.~Vinet, {\ita Rodrigues formulas for the 
Jack Polynomials and the Macdonald-Stanley conjecture
}  IMRN {\bf 9} 419-424 (1995).
\sas

\item{[10]}
L.~Lapointe and L.~Vinet, {\ita Rodrigues formulas for the {M}acdonald
  polynomials}, Adv. Math. {\bol 130} (1997), no.~2, 261--279.
\sas

\item {[11]}
I. G. Macdonald, {\ita  A new class of symmetric functions}, 
Actes du $20^e$ S\'eminaire Lotharingien, 
Publ. I.R.M.A. Strasbourg, (1988)
131-171.
\sas

\item {[12]}
I. G. Macdonald, {\ita Symmetric functions and Hall polynomials},
Second Edition, Clarendon Press, Oxford (1995).
\sas

\item{[13]}
S.~Sahi, {\ita Interpolation, integrality, and a generalization of {M}acdonald's
  polynomials}, Internat. Math. Res. Notices (1996), no.~10, 457--471.
\sas

\item {[14]} M. Zabrocki, {\ita q-Analogs of symmetric function operators},
(Preprint).

\end